\documentclass[11pt]{amsart}
\usepackage{amssymb}
\usepackage{amscd}
\usepackage[all]{xy}
\xyoption{ps}
\xyoption{dvips}
\usepackage{amsthm}

\newtheorem{Thm}{Theorem}[section]
\newtheorem{Lem}[Thm]{Lemma}
\newtheorem{Cor}[Thm]{Corollary}
\newtheorem{Prop}[Thm]{Proposition}

\newcommand{\md}{\operatorname{mod}}

\newcommand{\Hom}{\operatorname{Hom}}
\newcommand{\Ext}{\operatorname{Ext}}
\newcommand{\End}{\operatorname{End}}

\newcommand{\field}{k}
\newcommand{\spl}{{\rm Sp}}

\newcommand{\gldim}{\operatorname{gldim}}
\newcommand{\pdim}{\operatorname{projdim}}

\newcommand{\rdim}{\operatorname{repdim}}
\newcommand{\fdim}{\operatorname{findim}}
\newcommand{\cok}{\operatorname{Cok}}
\newcommand{\ke}{\operatorname{Ker}}
\newcommand{\im}{\operatorname{Im}}
\newcommand{\soc}{\operatorname{soc}}
\newcommand{\tp}{\operatorname{top}}
\newcommand{\add}{\operatorname{add}}

\begin{document}

\title{Radical embeddings and representation dimension}

\author{Karin Erdmann, Thorsten Holm, Osamu Iyama, Jan Schr\"oer}

\address{Karin Erdmann\newline
Mathematical Institute\newline
24-29 St. Giles\newline
Oxford OX1 3LB\newline 
ENGLAND}
\email{erdmann@maths.ox.ac.uk}

\address{Thorsten Holm\newline
Otto-von-Guericke-Universit\"at\newline
Institut f\"ur Algebra und Geometrie\newline
Postfach 4120\newline
39016 Magdeburg\newline
GERMANY}
\email{thorsten.holm@mathematik.uni-magdeburg.de}

\address{Osamu Iyama\newline
Department of Mathematics\newline
Himeji Institute of Technology\newline
Himeji, 671-2201\newline
JAPAN}
\email{iyama@sci.himeji-tech.ac.jp}

\address{Jan Schr\"oer\newline
Department of Pure Mathematics\newline
University of Leeds\newline
Leeds LS2 9JT\newline
ENGLAND}
\email{jschroer@maths.leeds.ac.uk}

\thanks{Mathematics Subject Classification (2000): 16E05, 16E10, 
16G20.}

\begin{abstract}
Given a representation-finite algebra $B$ and a 
subalgebra $A$ of $B$ such that the Jacobson radicals
of $A$ and $B$ coincide, we prove that the 
representation dimension of $A$ is at most three.
By a result of Igusa and Todorov, this implies that the finitistic
dimension of $A$ is finite.
\end{abstract}

\maketitle

\begin{center}
\vskip-0.3cm

DEDICATED TO THE MEMORY OF SHEILA BRENNER
\end{center}

\vskip0.5cm


\section{Introduction and main result}


Throughout, let $\field$ be a
field, and let $A$ be an associative finite-dimensional 
$\field$-algebra.
By $\md(A)$ we denote the category of finite-dimensional 
left $A$-modules.
A {\it generator-cogenerator} of $A$ is
an $A$-module $Z$ such that $A \oplus A^* \in \add(Z)$,
where $A^* = D(A_A)$ with $D = \Hom_\field(-,\field)$ the usual
duality.
Recall that $\add(Z)$ is the subcategory of all $A$-modules 
which are isomorphic to direct summands of finite direct sums of $Z$.
The {\it representation dimension} of $A$ is defined as
\[
\rdim(A) = \min \{ \gldim(\End_A(Z)) \mid Z \text{ a generator-cogenerator
of } A \}.
\]
Recall that an algebra $A$ is representation-finite
if there are only finitely many isomorphism classes of indecomposable 
$A$-modules.
Auslander \cite{A} showed that $A$ is representation-finite
if and
only if $\rdim(A) \leq 2$.

There are only few examples of algebras
where the representation dimension is known.
For an algebra $A$, we denote its Jacobson radical by $J_A$.
An algebra homomorphism $f: A \to B$ is called
a {\it radical embedding} if $f$ is a monomorphism
with $f(J_A) = J_B$.
Our main result is the following:

\vspace{0.2cm}
\begin{Thm}\label{mainthm}
If $f: A \to B$ is a radical embedding with $B$ a
representation-finite algebra, then $\rdim(A) \leq 3$.  
\end{Thm}
\vspace{0.2cm}

The {\it finitistic dimension} of $A$ is defined as
\[
\fdim(A) = \sup \{ \pdim(M) \mid M \in \md(A), \pdim(M) < \infty \}.
\]
It is a famous problem whether $\fdim(A) < \infty$
for all finite-dimensional algebras $A$.

If $A$ is an algebra with $\rdim(A) \leq 3$, then by a recent 
result of Igusa and Todorov \cite{IT} one gets that $\fdim(A) < \infty$.
This clearly underlines the importance of the representation dimension.
To check whether the finitistic dimension is finite is very difficult,
since (at least with the naive approach) we have to compute the projective
dimension of all modules.
To prove that the representation dimension of an algebra is at most 
three, one `just' has to guess a suitable module, which is a
generator-cogenerator and whose endomorphism algebra has global dimension
at most three.
Up to now there are no examples known where the representation dimension
is bigger than three.
It is shown in \cite{I_2} that $\rdim(A) < \infty$ for all algebras $A$.
For further results on the representation dimension of algebras
we refer to \cite{A}, \cite{X1} and \cite{X2}.
For proving that a particular algebra has
representation dimension at most three, the following result is 
often useful:

\vspace{0.2cm}
\begin{Prop}\label{projinj}
Let $A$ be a basic algebra, and let $P$ be an indecomposable 
projective-injective $A$-module. 
Define $B = A/\soc(P)$.
If $\rdim(B) \leq 3$, then $\rdim(A) \leq 3$.
\end{Prop}
\vspace{0.2cm}

For special biserial algebras, it was an open problem whether their
finitistic dimension is finite.
In \cite{S} it was proved that all special biserial algebras with 
at most two simple modules have finite finitistic dimension.
The following application of Theorem \ref{mainthm} settles this problem
for all special biserial algebras:

\vspace{0.2cm}
\begin{Cor}\label{maincor}
If $A$ is a special biserial algebra,
then we have $\rdim(A) \leq 3$ and $\fdim(A) < \infty$.
\end{Cor}
\vspace{0.2cm}

The proof of Theorem \ref{mainthm} yields an explicit construction
of a generator-cogenerator, whose endomorphism algebra has
the desired global dimension.
Namely, as before, let $f: A \to B$ be a radical embedding
with $B$ representation-finite.
Let $N$ be the direct sum of a complete set of representatives
of the isomorphism classes of indecomposable $B$-modules.
Note that we can consider $N$ also as an $A$-module.
Then define $C_f = A \oplus A^* \oplus N$.
We get the following result:

\vspace{0.2cm}
\begin{Thm}\label{thm2}
If $f: A \to B$ is a radical embedding with
$B$ a representation-finite algebra,
then $\End_A(C_f)$ is a quasi-hereditary algebra of
global dimension at most three.
\end{Thm}
\vspace{0.2cm}

The paper is organized as follows:
In Section \ref{proofs} we prove 
Theorem \ref{mainthm} and Proposition \ref{projinj}.
In Section \ref{rademb} we give a general construction principal
for radical embeddings.
This is applied in Section \ref{stringalg} 
to prove Corollary \ref{maincor}.
Section \ref{proofthm2} contains the proof of Theorem \ref{thm2}.
Finally, we discuss an example in Section \ref{example}.

In this paper, `modules' are finite-dimensional
left modules.
Although we often write maps on the left hand side,
we compose them as if they were on the right. 
Thus the composition of a map $\theta$ 
followed by a map $\phi$ is denoted $\theta\phi$.


\section{Proof of Theorem \ref{mainthm} and
Proposition \ref{projinj}}\label{proofs}


The proof of the following lemma is implicitly contained in
\cite[Chapter III, \S 3]{A}.
For convenience we repeat it here.

\vspace{0.2cm}
\begin{Lem}\label{lem1}
Let $A$ be an algebra, and let $M$ be a generator-cogenerator of
$A$.
Then for $n \geq 3$ the following are equivalent:
\begin{itemize}

\item[(1)]
For all indecomposable $A$-modules $X$ there
exists an exact sequence
\[
0 \to M_{n-2} \to \cdots \to M_1 \to M_0 \to X \to 0
\]
with $M_i \in \add(M)$, such that 
\begin{multline*}
0 \to \Hom_A(M,M_{n-2}) \to \cdots \to \Hom_A(M,M_1)\\
\to \Hom_A(M,M_0) \to \Hom_A(M,X) \to 0
\end{multline*}
is exact;

\item[(2)]
For all indecomposable $A$-modules $X$ there
exists an exact sequence
\[
0 \to X \to M_0 \to M_1 \to \cdots \to M_{n-2} \to 0 
\]
with $M_i \in \add(M)$, such that 
\begin{multline*}
0 \to \Hom_A(M_{n-2},M) \to \cdots \to \Hom_A(M_1,M)\\ 
\to \Hom_A(M_0,M) \to \Hom_A(X,M) \to 0
\end{multline*}
is exact;

\item[(3)]
$\gldim(\End_A(M)) \leq n$.

\end{itemize}
\end{Lem}

\begin{proof}
For brevity set $E = \End_A(M)$.
Assume that (1) holds.
Let $T$ be an $E$-module, and let
\[
\Hom_A(M,M'') \stackrel{F}{\longrightarrow} \Hom_A(M,M') \to T \to 0
\]
be a projective presentation of $T$.
Then $F = \Hom_A(M,f)$ for some homomorphism $f: M'' \to M'$.
Thus we get an exact sequence
\[
0 \to \ke(f) \to M'' \to M'.
\]
Using our assumption, we get an exact sequence
\[
0 \to M_{n-2} \to \cdots \to M_1 \to M_0 \to \ke(f) \to 0
\]
having the properties described in (1), here we set $X = \ke(f)$.
This yields an exact sequence
\[
0 \to M_{n-2} \to \cdots \to M_1 \to M_0 \to M'' \to M'.
\]
Applying $\Hom_A(M,-)$ gives an exact sequence
\begin{multline*}
0 \to \Hom_A(M,M_{n-2}) \to \cdots \to \Hom_A(M,M_1) \to \Hom_A(M,M_0)\\
\to \Hom_A(M,M'') \to \Hom_A(M,M') \to T \to 0.
\end{multline*}
Thus $\pdim(T) \leq n$ for all $E$-modules $T$.
We get $\gldim(E) \leq n$.
Thus (3) is true.

Next, assume that (3) holds.
For an $A$-module $X$, let
\[
0 \to X \to I_0 \stackrel{h}{\longrightarrow} I_1
\]
be an injective presentation.
Note that $I_0,I_1 \in \add(M)$.
We get an exact sequence
\[
0 \to \Hom_A(M,X) \to \Hom_A(M,I_0) \to \Hom_A(M,I_1) \to Y \to 0
\] 
with $Y = \cok(\Hom_A(M,h))$.
Now $\Hom_A(M,X)$ is the second syzygy module $\Omega^2(Y)$ of
$Y$.
Since $\gldim(E) \leq n$, we know that $\pdim(\Omega^2(Y)) \leq n-2$.
Thus there exists an exact sequence
\begin{multline*}
0 \to \Hom_A(M,M_{n-2}) \to \cdots \to \Hom_A(M,M_1)\\
\to \Hom_A(M,M_0) \to \Hom_A(M,X) \to 0
\end{multline*}
with $M_i \in \add(M)$.
This yields an exact sequence
\[
0 \to M_{n-2} \to \cdots \to M_1 \to M_0 \to X \to 0.
\]
Thus (1) follows.
The equivalence of the statements (2) and (3) is proved dually.
This finishes the proof.
\end{proof}

Let $B$ be an algebra, and let $A \subseteq B$ be a subalgebra of $B$.
We regard any $B$-module also as an $A$-module in the obvious way.
For an $A$-module $X$, define $X^- = \Hom_A(B,X)$.
Furthermore, we identify $X$ and $\Hom_A(A,X)$.
Let 
\[
\epsilon_X: X^- \to X 
\] 
be the natural map induced by the inclusion $A \subseteq B$.
Note that $\epsilon_X$ is an $A$-module homomorphism.
Observe also that $X^-$ is a $B$-module.

\vspace{0.2cm}
\begin{Lem}\label{lem2}
Let $A$ be a subalgebra of an algebra $B$, and
let $X$ be an $A$-module.
Then 
\[
\Hom_B(Y,X^-) \to \Hom_A(Y,X),
\]
\[
f \mapsto f \epsilon_X
\]
is an isomorphism for all $B$-modules $Y$.
\end{Lem}

\begin{proof}
For all $B$-modules $Y$ we have
\[
\Hom_B(Y,\Hom_A(B,X)) \cong \Hom_A(B \otimes_B Y,X) \cong \Hom_A(Y,X),
\]
where the isomorphisms are given by
\[
f \mapsto (b \otimes y \mapsto f(y)(b)) \mapsto (y\mapsto f(y)(1)).
\]
Thus the composition maps $f$ to $f\epsilon_X$.
\end{proof}

\begin{Lem}\label{lem3}
Let $A$ be a subalgebra of an algebra $B$
such that $J_A = J_B$.
Then $\cok(\epsilon_X)$ and $\ke(\epsilon_X)$ are
semisimple $A$-modules for all $A$-modules $X$.
\end{Lem}

\begin{proof}
The sequence
\[
\xymatrix{
0 \ar[r] & \Hom_A(B/A,X) \ar[r] & X^- \ar[r]^{\epsilon_X} & X 
\ar[r] & \Ext^1_A(B/A,X) 
}
\]
is exact.
We have $(B/A)J_A = 0$.
Thus we get $J_A \Ext^1_A(B/A,X) = 0$, which
implies that $\Ext^1_A(B/A,X)$ is a semisimple $A$-module.
Thus $\cok(\epsilon_X)$ is a semisimple $A$-module, since
it is a submodule of $\Ext^1_A(B/A,X)$.
Also $\ke(\epsilon_X)$ is semisimple, since
$J_A \Hom_A(B/A,X) = 0$.
\end{proof}

\begin{proof}[Proof of Theorem \ref{mainthm}]
Assume that $B$ is a representation-finite algebra,
and let $M_1, \cdots, M_n$ be a complete set of representatives of 
isomorphism classes of indecomposable $B$-modules.
Without loss of generality we assume that $A$ is a subalgebra of
$B$ such that $J_A = J_B$.
Define $N = \bigoplus_{i=1}^n M_i$, and let
$M = A \oplus A^* \oplus N$.

We claim that $\gldim(\End_A(M)) \leq 3$.
To prove this, we use the criterion presented in Lemma \ref{lem1}.

If $X$ is an indecomposable injective $A$-module, then 
we get a short exact sequence 
\[
0 \to 0 \to X \to X \to 0.
\]
Setting $M_0 = X$ and $M_1 = 0$, we see that this sequence
satisfies the conditions in Lemma \ref{lem1}(1).

Assume next that $X$ is an indecomposable non-injective $A$-module.
We know by Lemma \ref{lem3} that $\cok(\epsilon_X)$ is a semisimple
$A$-module.
By $\pi: P \to \cok(\epsilon_X)$ we denote the projective cover 
of $\cok(\epsilon_X)$.
Since $P$ is a projective $A$-module, there
exists a homomorphism $p: P \to X$ such that the diagram
\[
\xymatrix{
 & 0 \ar[r] & J_A P \ar[r]^\iota & P \ar[r]^-{\pi} \ar[d]^p & 
\cok(\epsilon_X) \ar[r] \ar@{=}[d]& 0\\
0 \ar[r] & \ke(\epsilon_X) \ar[r] & X^- \ar[r]^{\epsilon_X} & X 
\ar[r] & \cok(\epsilon_X) \ar[r] & 0 
}
\]
commutes and has exact rows.
Observe that the map
\[
\bigl( \begin{smallmatrix}
\epsilon_X\\
p
\end{smallmatrix} \bigr): X^- \oplus P \to X 
\]
is an epimorphism of $A$-modules.
Note also that $X^- \oplus P \in \add(M)$.
We take
\[
M_0 = X^- \oplus P
\]
and will show that this works.
It follows from Lemma \ref{lem2} that the map
\[
\Hom_A(Y,X^- \oplus P) \to \Hom_A(Y,X)
\]
\[
(f,g) \mapsto (f \epsilon_X + g p) 
\]
is surjective for any $B$-module $Y$. 
Since $A$ is projective, it follows that the map
\[
\Hom_A(A,X^- \oplus P) \to \Hom_A(A,X)
\]
\[
(f,g) \mapsto (f \epsilon_X + g p) 
\]
is surjective as well.

Finally, take an injective $A$-module $I$, and some
homomorphism $f \in \Hom_A(I,X)$.
Let $I \to I/\soc(I)$ be the canonical projection.
Since $X$ is not injective, 
there exists a homomorphism $g: I/\soc(I) \to X$ such that
the diagram 
\[
\xymatrix{
I \ar[r] \ar[d]_f & I/\soc(I) \ar[dl]^g\\
X & 
}
\]
commutes.
We have $I = \bigoplus_{i=1}^t (e_iA)^*$ for some
primitive idempotents $e_i$ in $A$.
Since $A \subseteq B$, the $e_i$ are also idempotents of $B$.
From $J_A = J_B$ we get 
$I/\soc(I) = \bigoplus_{i=1}^t (e_iJ_B)^*$.
Thus $I/\soc(I)$ is also a $B$-module.
Thus by Lemma \ref{lem2}, $g$ factors through $\epsilon_X$.

Altogether, we proved that for any $A$-module $Z \in \add(M)$ and
any homomorphism $f: Z \to X$ of $A$-modules with $X$ indecomposable 
there exists a homomorphism $g: Z \to X^- \oplus P$ of $A$-modules, 
such that the diagram
\[
\xymatrix{
& Z \ar[dl]_g \ar[d]^f \\
X^- \oplus P \ar[r]_-{ \bigl( 
\begin{smallmatrix}
\epsilon_X\\
p
\end{smallmatrix} \bigr)} & X
}
\]
commutes.
Next, we show that the kernel of the map
$\bigl( \begin{smallmatrix}
\epsilon_X\\
p
\end{smallmatrix} \bigr): M_0 \to X$ 
belongs to $\add(M)$.

Let $\epsilon_X': X^- \to \im(\epsilon_X)$ be the 
epimorphism induced by $\epsilon_X$.
There are the obvious inclusion maps $\im(\epsilon_X) \to X$ and
$J_A P \to P$.
Clearly, there exists a homomorphism $p': J_A P \to \im(\epsilon_X)$
such that the diagram
\[
\xymatrix{
J_A P \ar[r]^\iota \ar[d]^{p'} & P \ar[d]^p\\
\im(\epsilon_X) \ar[r] & X
}
\]
commutes.
Now we construct the pullback of $p'$ and get the commutative diagram
\[
\xymatrix{
0 \ar[r] & \ke(\epsilon_X) \ar[r] \ar@{=}[d] & Y \ar[r] \ar[d] & J_A P 
\ar[d]^{p'} \ar[r] & 0\\
0 \ar[r] & \ke(\epsilon_X) \ar[r] & X^- \ar[r]^-{\epsilon_X'} 
& \im(\epsilon_X) \ar[r] & 0.
}
\]
Thus
\[
Y = \ke \bigl( \begin{smallmatrix}
\epsilon_X\\
p
\end{smallmatrix} \bigr).
\]
We have $P = \bigoplus_{i=1}^t Ae_i$ for  some
primitive idempotents $e_i$ in $A$.
Since $J_A = J_B$, we get
$J_AP = \bigoplus_{i=1}^t J_Be_i$. 
Thus $J_AP$ is also a $B$-module.
Thus, by Lemma \ref{lem2}, the homomorphism $p'$ factors through
$\epsilon_X'$.
Thus the short exact sequence 
\[
0 \to \ke(\epsilon_X) \to Y \to J_A P \to 0
\]
splits, and we get
$Y \cong \ke(\epsilon_X) \oplus J_A P$.
By the construction of $M$ this implies $Y \in \add(M)$.
Now set $M_1 = Y$.

Thus for each $A$-module $X$ we constructed a short exact sequence
\[
0 \to M_1 \to M_0 \to X \to 0
\]
with the properties required in Lemma \ref{lem1}.
We get $\gldim(\End_A(M)) \leq 3$.
Since $M$ is a generator-cogenerator of $A$, it follows that 
$\rdim(A) \leq 3$.
This finishes the proof.
\end{proof}

\noindent
{\em Remark.} Theorem 1.1 and its proof hold 
under the weaker assumption that $f$ is a monomorphism such that
$f(J_A)$ is a two-sided ideal of $B$ (not necessarily equal to $J_B$). 
So far we are not aware of interesting applications of this
slightly more general result. Thus we refrain from giving details here.

\begin{proof}[Proof of Proposition \ref{projinj}]
Next, we prove Proposition \ref{projinj}.
We have $B = A/\soc(P)$.
Thus there is a surjective algebra homomorphism
$f: A \to B$.
We can regard any $B$-module as an $A$-module with the $A$-module
structure induced by $f$.

Let $N$ be a generator-cogenerator of $B$ with
$\gldim(\End_B(N)) \leq 3$.
Define $M = N \oplus P$.
Observe that $M$ is a generator-cogenerator of $A$.
We claim that $\gldim(\End_A(M)) \leq 3$.
To check this, we use again Lemma \ref{lem1}.
Let $X$ be any indecomposable $A$-module.
If $X = P$, then we get a short exact sequence 
\[
0 \to 0 \to X \to X \to 0.
\] 
Set $M_0 = X$ and $M_1 = 0$.
It is easy to verify that this sequence
satisfies the conditions required in Lemma \ref{lem1}.
Next, assume that $X$ is not isomorphic to $P$.
Thus $X$ is an indecomposable $B$-module.
Applying Lemma \ref{lem1} and our assumption 
$\gldim(\End_B(N)) \leq 3$, we get a short exact sequence
\[
0 \to N_1 \to N_0 \to X \to 0
\]
of $B$-modules with $N_0,N_1 \in \add(N)$ and
\[
0 \to \Hom_B(N,N_1) \to \Hom_B(N,N_0) \to \Hom_B(N,X) \to 0
\]
an exact sequence.
Since $P$ is projective, the functor $\Hom_A(P,-)$ is exact.
Thus we get an exact sequence
\[
0 \to \Hom_A(M,N_1) \to \Hom_A(M,N_0) \to \Hom_A(M,X) \to 0.
\]
This enables us to apply Lemma \ref{lem1} again, and we get
$\gldim(\End_A(M)) \leq 3$.
This finishes the proof.
\end{proof}


\section{Construction of radical embeddings}\label{rademb}


A {\it quiver} is a quadruple $Q = (Q_0,Q_1,s,e)$,
where $Q_0$ and $Q_1$ are finite sets and 
$s,e: Q_1 \to Q_0$ are maps.
We call the elements in $Q_0$ the {\it vertices} of $Q$,
and the elements in $Q_1$ the {\it arrows} of $Q$.
A {\it path} of length $n \geq 1$ in $Q$ is of the form 
$\alpha_1 \alpha_2 \cdots \alpha_n$ where the $\alpha_i$ are arrows with 
$s(\alpha_i) = e(\alpha_{i+1})$ for $1 \leq i \leq n-1$. 
Additionally, there is a path $e_i$ of length zero for each vertex 
$i \in Q_0$.
By $\field Q$ we denote the path algebra of $Q$ with basis the
set of all paths in $Q$.
The multiplication is given by concatenation of paths.

A {\it relation} for $Q$ is a linear combination
$\sum_{i=1}^t \lambda_i r_i$ such that 
$\lambda_i \in \field^*$ and
the $r_i$ are paths of the form $\alpha_i p_i \beta_i$
with $\alpha_i,\beta_i \in Q_1$ 
such that $s(\beta_i) = s(\beta_j)$ and $e(\alpha_i) = e(\alpha_j)$ for all 
$1 \leq i,j \leq t$.

A {\it basic algebra} is a finite-dimensional algebra of the form
$\field Q/I$, where the ideal $I$ is generated by a set of relations.
By a result of Gabriel, any finite-dimensional $\field$-algebra is
Morita equivalent to a basic algebra provided we assume that $\field$
is algebraically closed.

Now, let $A = \field Q/I$ be a basic algebra
with $Q = (Q_0,Q_1,s,e)$.
Let $l \in Q_0$ be a vertex.
Define 
\[
S(l) = \{ \alpha \in Q_1 \mid s(\alpha) = l \}
\]
and 
\[
E(l) = \{ \beta \in Q_1 \mid e(\beta) = l \}.
\]
Note that the intersection of $S(l)$ and $E(l)$
might be non-empty.

Let $S(l) = S_1 \cup S_2$ and $E(l) = E_1 \cup E_2$
be disjoint unions.
We call $(S_1,S_2,E_1,E_2)$ a {\it splitting datum at} $l$ if 
the following hold:
\begin{itemize}

\item[(1)]
For $\alpha \in S_i$ and $\beta \in E_j$ we have
$\alpha \beta = 0$ whenever $i \not= j$;

\item[(2)]
The ideal $I$ can be generated by a set $\rho$ of relations
of the form $\sum_{i=1}^t \lambda_i \alpha_ip_i\beta_i$
such that
$\{ \alpha_i \mid 1 \leq i \leq t \} \cap E_j = \emptyset$
for $j = 1$ or $j=2$, 
and 
$\{ \beta_i \mid 1 \leq i \leq t \} \cap S_j = \emptyset$ 
for $j=1$ or $j=2$.

\end{itemize}
Note that condition (2) in the above definition is automatically
satisfied, if we assume that $I$ is a monomial ideal, i.e. if 
$I$ can be generated by a set of paths in $Q$.
Given a splitting datum $\spl = (S_1,S_2,E_1,E_2)$ at $l$, 
we construct from $Q$ a new quiver 
\[
Q^\spl = (Q_0^\spl,Q_1^\spl,s^\spl, e^\spl)
\]
as follows:
Let 
\[
Q_0^\spl = \{ l_1, l_2 \} \cup Q_0 \setminus \{ l \},
\]
and set $Q_1^\spl = Q_1$.
The maps $s^\spl, e^\spl: Q_1^\spl \to Q_0^\spl$
are 
\[
s^\spl(\alpha) = 
\begin{cases}
s(\alpha) & \text{if $s(\alpha) \not= l$},\\
l_1       & \text{if $\alpha \in S_1$},\\
l_2       & \text{if $\alpha \in S_2$},
\end{cases}
\]
and 
\[
e^\spl(\alpha) = 
\begin{cases}
e(\alpha) & \text{if $e(\alpha) \not= l$},\\
l_1       & \text{if $\alpha \in E_1$},\\
l_2       & \text{if $\alpha \in E_2$}.
\end{cases}
\]
Let $\rho$ be a set of relations for $Q$ which satisfy 
condition (2) above.
Define 
\[
\rho^\spl = \rho \setminus \{ \alpha\beta \mid \alpha \in S_i,
\beta \in E_j, i \not= j \}. 
\]
Then each element in $\rho^\spl$ is also a relation for the quiver
$Q^\spl$.
Let $I^\spl$ be the ideal of $\field Q^\spl$
generated by the relations in $\rho^\spl$.
Set 
\[
A^\spl = \field Q^\spl/I^\spl.
\]
We get the following result:

\vspace{0.2cm}
\begin{Lem}\label{splitting}
Let $A = \field Q/I$ be a basic algebra, and let
$\spl$ be a splitting datum at some vertex of $Q$.
Then there exists a radical embedding
\[
A \to A^\spl.
\]
\end{Lem}

\begin{proof}
Let $\spl$ be a splitting datum at some vertex $l \in Q_0$.
We construct a map $f: A \to A^\spl$ as follows:
For the arrows $\alpha \in Q_1$ we just define
\[
f(\alpha) = \alpha.
\]
For a vertex $j \in Q_0$ let
\[
f(e_j) = 
\begin{cases}
e_j & \text{if $j \not= l$},\\
e_{l_1} + e_{l_2} & \text{if $j = l$}.
\end{cases}
\]
It follows directly from the definition of a splitting
datum that $f$ can be extended to an algebra homomorphism.
It is also clear that $f$ is a monomorphism and satisfies
$f(J_A) = J_{A^\spl}$.
\end{proof}

The above lemma is useful for the construction of radical
embeddings.
In fact, in can be applied to numerous situations.
In the next section, we illustrate this for one of the
most important classes of tame algebras, the string algebras.


\section{Proof of Corollary \ref{maincor}}\label{stringalg}


A basic algebra $A = kQ/I$ is called a {\it special biserial algebra}
if the following hold:
\begin{enumerate}

\item[(1)] Any vertex of $Q$ is the starting point of at most two 
arrows and also the end point of at most two arrows;

\item [(2)] Let $\beta$ be an arrow in $Q_1$.
Then there is at most one arrow 
$\alpha$ with $\alpha \beta \notin I$ 
and at at most one arrow $\gamma$ with  
$\beta \gamma \notin I$;

\item[(3)] There exists some $N$ such that each path of length
at least $N$ lies in $I$, i.e. $A$ is finite-dimensional.

\end{enumerate}
A {\it string algebra} is a special biserial algebra $\field Q/I$ 
which satisfies the additional condition that $I$ is generated
by paths.
For details and further references on string algebras we refer
to \cite{BR}.

\begin{proof}[Proof of Corollary \ref{maincor}]
Let $A = kQ/I$ be a string algebra.
Define
\[
c(A) = |\{ l \in Q_0 \mid |S(l)| = 2| + |\{ l \in Q_0 \mid |E(l)| = 2 \}|.
\]
If $c(A) = 0$,
then $Q$ is a disjoint union of quivers which are of type
$\mathbb{A}$ with linear orientation or of type 
$\widetilde{\mathbb{A}}$ with cyclic orientation.
But string algebras with such underlying quivers
are representation-finite.
In fact, it is easy to check that for a string algebra
$A$ all indecomposable $A$-modules are serial if and only if $c(A) = 0$.

Thus, assume $c(A) \geq 1$.
Let $l \in Q_0$ such that $|S(l)| = 2$ or $|E(l)| = 2$.
First, we consider the case $|S(l)| = 2$, say
$S(l) = \{ \alpha_1, \alpha_2 \}$.
We construct a splitting datum $\spl = (S_1,S_2,E_1,E_2)$ at $l$ 
as follows:
Let $S_1 = \{ \alpha_1 \}$, $S_2 = \{ \alpha_2 \}$,
$E_1 = \{ \beta \in E(l) \mid \alpha_2 \beta = 0 \}$ 
and
$E_2 = E(l) \setminus E_1$.
It follows directly from the definition of a string algebra,
that $\spl$ is a splitting datum.
Now $A^\spl$ is again a string algebra, and we have
\[
c(A^\spl) \leq c(A) - 1.
\]
The case $|E(l)| = 2$ is done in the same way.
Repeating this construction a finite number of times and applying
Lemma \ref{splitting} yields
a chain 
\[
A = A_1 \to A_2 \to \cdots \to A_t = B
\]
of radical embeddings, where $B$ is a string
algebra with $c(B) = 0$.
As observed above, $B$ is representation-finite.
Thus, for any string algebra $A$, we get a radical embedding
$A \to B$ with $B$ representation-finite.
Then Theorem \ref{mainthm} yields that $\rdim(A) \leq 3$.

Next, assume that $A$ is a special biserial algebra.
Then we get from $A$ to a string algebra $B$ by successively factoring out
socles of indecomposable projective-injective modules.
Applying Proposition \ref{projinj} after each step, we get
$\rdim(A) \leq 3$.
Now we use the result in \cite{IT} and get $\fdim(A) < \infty$
for any special biserial algebra $A$.
This finishes the proof.
\end{proof}

Note that for $A$ a string algebra,
the proofs of Theorem \ref{mainthm} and Corollary \ref{maincor}
yield an explicit construction of a generator-cogenerator
$M$ of $A$ such that $\gldim(\End_A(M)) \leq 3$.
Namely, take $M$ as the direct sum of a
complete set of representatives of isomorphism classes of string modules,
which are projective, injective or serial.


\section{Proof of Theorem \ref{thm2}}\label{proofthm2}


Let $A$ be a subalgebra of an algebra $B$.
We have the `induction' functor
\[
T = {_BB} \otimes_A -: \md(A) \to \md(B),
\]
which is left adjoint to the `inclusion' functor
\[
F = \Hom_B(B, -): \md(B) \to \md(A).
\]
Thus for any $A$-module $Y$ and any $B$-module $X$ 
we get an isomorphism
\[
\phi_{X,Y}: \Hom_B(TY,X) \to \Hom_A(Y,FX).
\]
For the sake of brevity we will just write $\phi$ instead if $\phi_{X,Y}$.
Sometimes we will omit writing $F$. 
Let
\[
e: F T \to 1_{\md(B)}
\]
be the corresponding counit, so that
\[
e_X = \phi^{-1}(1_{FX}): B \otimes_A \Hom_B(B,X) = T(FX) \to X
\]
is just the multiplication map.
This is a $B$-homomorphism.
The unit of this adjunction is the natural transformation
\[
\delta: 1_{\md(A)} \to TF,
\]
so that for $Y \in \md(A)$ we have 
\[
\delta_Y = \phi(1_{TY}): Y \to  F(TY) 
\]
\[
y \mapsto (1 \otimes y)
\]
if $F(TY)=\Hom_B(B, B\otimes_AY)$ is identified with $B\otimes_AY$.
This is an $A$-module homomorphism.
Note also that $\phi^{-1}(g) = T(g) e_X$ for
$g: Y \to FX$ an $A$-homomorphism, and
$\phi(f) = \delta_Y F(f)$ for $f: TY \to X$ a $B$-homomorphism.

\vspace{0.2cm}
\begin{Lem}\label{lemma5}
Let $A$ be a subalgebra of an algebra $B$ such that $J_A = J_B$.
If $X$ is a $B$-module, then
as a $B$-module, we have $T(FX) \cong B\otimes_AFX\cong X \oplus S$ where
$S$ is a semisimple $B$-module.
\end{Lem}

\begin{proof}
Write $Y = FX$. 

(1) First consider the exact sequences and the resulting commutative diagram 
\[
\xymatrix{
0 \ar[r] & Y \ar[r]^-{\delta_Y}\ar[d] & B\otimes_A Y \ar[r]\ar[d] & B/A
\otimes_A Y \ar[r]\ar@{=}[d] & 0 \\
( 0 \ar[r] &) \tp(Y)  \ar[r] & \tp(B\otimes_AY) \ar[r]  &
\tp(B/A\otimes_AY)\ar[r] & 0
}
\]
of $A$-homomorphism, obtained by taking radical quotients
(here $\tp(M) = M/J_AM$), where 
the vertical maps are the canonical epimorphisms.
Since $\tp(Y)$ is the restriction of a $B$-module, the map
$Y \to \tp(Y)$ factors through $\delta_Y$, see the dual of Lemma \ref{lem2},
that is by adjointness.
Hence the lower row is a split short exact sequence.

(2)
Let $e_X$ be the counit of the adjunction.
Then we have a commutative diagram
\[
\xymatrix{
0 \ar[r] & \ke(e_X)  \ar[r] \ar[d]^{p} & B \otimes_A Y \ar[r]^{e_X} \ar[d] &
X\ar[r]  \ar[d]& 0\\
( 0 \ar[r] &) \tp(\ke(e_X)) \ar[r]  & \tp(B \otimes_A Y) \ar[r] & 
\tp(X) \ar[r]& 0
}
\]
of $B$-homomorphisms, which has exact rows.
By $l(M)$ we denote the length of a module $M$.
We have (using that the lower sequence in (1) is split exact)
\begin{multline*}
l(\ke(e_X)) = l(B\otimes_AY) - l(X) = 
l(B \otimes_A Y) - l(Y)\\
= l(B/A \otimes_A Y) = l(\tp(B \otimes_A Y)) - l(\tp(Y))\\
= l(\tp(B\otimes_AY)) - l(\tp(X)) \leq l(\tp(\ke(e_X))).
\end{multline*}
Thus $\ke(e_X)$ is a semisimple $B$-module, and
$p$ is an isomorphism, and both rows in the diagram in (2) are split short 
exact sequences of $B$-modules.
\end{proof}

\begin{proof}[Proof of Theorem \ref{thm2}]
Let $f: A \to B$ be a radical embedding with $B$ a representation-finite
algebra.
We assume without loss of generality that $f$ is an inclusion map,
i.e. $A \subseteq B$ with $J_A = J_B$.
Let $N$ be the direct sum of a complete set of representatives
of the isomorphism classes of indecomposable $B$-modules, say
$N = \bigoplus_i N_i$ with $N_i$ indecomposable.
Define $M = A \oplus A^* \oplus FN$.
From the proof of Theorem \ref{mainthm} we already know that
$\gldim(\End_A(M)) \leq 3$.
We claim that $\End_A(M)$ is quasi-hereditary.

(1)
Set $\Gamma = \End_A(M)$.
Recall that the isomorphism classes of simple modules of
the endomorphism algebra of a module are indexed
by the isomorphism classes of its indecomposable direct summands.
Let 
\[
R = \End_B(B \otimes_A M) = \End_B(TM) = \End_B(TA \oplus TA^* \oplus T(FN)).
\]
Lemma \ref{lemma5} implies
$\add(B \otimes_A M) = \add(N)$.
Since $B$ is a representation-finite algebra, it follows from
\cite[Chapter III, \S 4]{A} that $\gldim(R) \leq 2$.
This implies that $R$ is a quasi-hereditary algebra with respect to some
partial order $\leq_R$ so that the labels given by the simple direct
summands of $N$ are maximal, see \cite{DR}.

Define a partial order $\leq$ on the labels for the simple
$\Gamma$-modules as follows:
Let $X$ and $Y$ be non-isomorphic indecomposable direct summands of the 
$A$-module $M$.
Set $X < Y$ if and only if one of the following holds:
\begin{itemize}

\item $X \cong FN_i$ and $Y \cong FN_j$ for some $i,j$ such that
$N_i <_R N_j$;

\item $X$ is not isomorphic a direct summand of $FN$, and 
$Y \cong FN_i$ for some $i$.

\end{itemize}
Note that this is a partial order:
The only indecomposable $B$-modules, which
could become isomorphic as $A$-modules, are simple ones.
Namely, if $N_i$ and $N_j$ are $B$-modules with $FN_i \cong FN_j$,
then $T(FN_i) \cong T(FN_j)$, 
and if they are not simple, then Lemma \ref{lemma5} implies that  
$N_i$ and $N_j$ are isomorphic.
It follows that all simple direct summands of $M$ are maximal with
respect to $\leq$.

For any indecomposable direct summand $X$ of $M$, we have
the submodule $U(X)$ of the projective $\Gamma$-module 
$P(X) = \Hom_A(M, X)$, which is defined to be
the span of all homomorphisms $M \to X$, which factor through some
$Y$ with $Y > X$.
The {\it standard module} associated to $X$ is defined
as $\Delta(X) = P(X)/U(X)$.
By $L(X)$ we denote the top of $P(X)$.
Thus $L(X)$ is simple.

We have to show that for each $X$, the module $P(X)$ has a filtration by
standard modules with $\Delta(X)$ occurring only once, 
and if $\Delta(Y)$ occurs, then $Y \geq X$.

(2)
For $X$ simple we have $\Delta(X) = P(X)$.
Assume now that $X$ is not isomorphic to a direct summand of $FN$.
Thus $X$ is a projective or injective $A$-module (and not simple).
In case $X$ is projective, the radical of $X$ is of the form
$\bigoplus_i X_i$ where $X_i = FX_i'$ with $X_i'$ an indecomposable
$B$-module for all $i$.
Then we get the exact sequence
\[
0 \to \bigoplus_i P(X_i) \to P(X),
\]
and the cokernel is one-dimensional, hence is $L(X)$.
Since $X_i >  X$ it follows that $\Delta(X) = L(X)$.

In the second case, we have a short exact sequence
of $A$-modules
\[
0 \to \Hom_A(B/A, X) \to X^-  \stackrel{\epsilon_X}{\longrightarrow} X \to 0
\]
with the kernel a  semisimple $A$-module,
write it as $\bigoplus_i S_i$ with $S_i$ simple.
Here we use Lemma \ref{lem3}.
Write also $X^- = \bigoplus_j X_j$ where $X_j = FX_j'$ 
with $X_j'$ an indecomposable $B$-module for all $j$.
This gives an exact sequence
\[
0 \to \bigoplus_i \Delta(S_i) \to \bigoplus_j P(X_j)
\to P(X).
\]
We claim that the cokernel at $P(X)$ is simple.
Suppose $g: W \to X$ is an $A$-homomorphism where $W$ is an indecomposable
direct summand of $M$.
If $W$ is isomorphic to a direct summand of $FN$, then
$g$ factors through $\epsilon_X$, via adjointness, and clearly it 
factors if $W$ is projective.
Suppose $W$ is indecomposable injective and not isomorphic to a 
direct summand of $FN$.
If $g$ is not an isomorphism, then it factors through the socle 
quotient of $W$.
But this is of the form $FW'$ for some $B$-module $W'$. 
Hence $g$ factors through $\epsilon_X$ again.
It follows that the cokernel is
$L(X)$ and is isomorphic to $\Delta(X)$.

(3)
So assume now that $X = FX'$ with $X'$ an indecomposable $B$-module,
which is not simple.
Let $\phi$ be the adjoint isomorphism
\[
\phi: \Hom_B(TM,X') \cong \Hom_A(M,X) = P(X). 
\]
Through the ring homomorphism
$T: \Gamma \to R$,
$\phi^{-1}$ induces $\Gamma$-isomorphisms
$P(X) \to P_{R}(X)$,
$U(X) \to U_{R}(X)$ and
$\Delta(X) \to \Delta_{R}(X)$.

We only have to show that this is compatible with factorizing through
some module $Z \in \add(N)$, modulo maps factorizing through
a semisimple module.

Suppose that $g: M \to X$ is an $A$-homomorphism 
with a factorization $g = \alpha \beta$, where $\alpha: M \to FZ$
and $\beta: FZ \to X$ are $A$-homomorphisms. 
Then we have
\[
\phi^{-1}(g) = T(\alpha)\phi^{-1}(\beta): TM \to X'.
\]
Hence $\phi^{-1}(g)$ factors through $T(FZ)$.
Since $Z$ is a $B$-module, Lemma \ref{lemma5} implies 
that $T(FZ) \cong Z \oplus S$
as a $B$-module with $S$ semisimple, and this is what we need.

Conversely, suppose $f: TM \to X'$ is a $B$-module homomorphism,
which factors through a $B$-module $Z$, say $f = \alpha\beta$, where
$\alpha: TM \to Z$ and $\beta: Z \to X'$. 
Then $\phi(f) = \phi(\alpha)F(\beta)$,  hence it
factors through $FZ$.

Since $R$ is quasi-hereditary with respect to $\leq_R$, each
indecomposable projective module $\Hom_B(TM,X')$ has
a filtration by standard modules of the right kind.
The above shows that the adjoint isomorphism identifies this 
filtration with a filtration of $\Hom_A(M, X)$ by standard modules 
for $\Gamma$.

It remains to show that  $L(X)$ occurs with multiplicity one as
a composition factor of $\Delta(X)$.
This is clear if $X$ is simple, and we have already
seen it in case $X$ is not isomorphic to a direct summand of $FN$.
If $X$ is isomorphic to a direct summand of $FN$,
then this multiplicity is the same as the multiplicity of
$L_{R}(X)$ in $\Delta_{R}(X)$, hence it is one.
\end{proof}


\section{An Example}\label{example}


Let $A = \field Q/I$ where $Q$ is the quiver
with one vertex $x$ and two loops $a, b$ and
$I = (a^2, b^2, (ab)^2, (ba)^2 )$.
When the field has characteristic 2 this is the socle quotient of
the group algebra of the dihedral group of order $8$.

Let $\spl = (S_1,S_2,E_1,E_2)$ where
$S_1 = \{ a \}$, $S_2 = \{ b \}$, $E_1 = \{ b \}$ and $E_2 = \{ a \}$.
Clearly, $\spl$ is a splitting datum at $x$.
Following the general construction in Section \ref{rademb}, 
we get the quiver $Q^\spl$ with two vertices $l_1$ 
and $l_2$ and arrows $a: l_1 \to l_2$ and $b: l_2 \to l_1$.
The ideal $I^\spl$ is generated by all paths of length 4 in $Q^\spl$.
The algebra $A^\spl$ is a Nakayama algebra.
Every indecomposable $A^\spl$-module is serial, and visibly its restriction 
to $A$ remains serial.

Hence $M = A \oplus A^* \oplus N$, where $N$ is the direct sum
of a complete set of representatives of isomorphism classes of 
serial string modules over $A$.
We denote these string modules as $M(C)$ for $C$ in 
$\{ 1_x, a, b, ab, ba, aba, bab \}$,
and we write $\field = M(1_x)$ for the simple $A$-module.
(For example, $M(ab)$ has basis $\{ v, bv, abv \}$).

Let $\Gamma = \End_A(M)$. 
We can see directly that $\gldim(\Gamma) = 3$:
For each indecomposable direct summand $W$ of $M$, we write $P(W)$ for the 
indecomposable projective $\Gamma$-module $\Hom_A(M,W)$. 
Let $L(W)$ be the simple top of $P(W)$.

(1) The radical $J_A$ belongs to $\add(M)$, and 
the inclusion $J_A \to A$
gives an inclusion of projective $\Gamma$-modules 
\[
0 \to \Hom_A(M, J_A) \to \Hom_A(M, A) = P(A).
\]
Clearly, the cokernel is 1-dimensional, hence it is the simple module
$L(A)$.
This implies $\pdim(L(A)) \leq 1$.

(2) We have an exact sequence
\[
0 \to \field \to M(aba) \oplus M(bab) \to A^*  \to 0.
\]
This gives an exact sequence 
\[
0 \to P(\field) \to P(M(aba)) \oplus P(M(bab)) \to P(A^*).
\]
of $\Gamma$-modules.
We claim that the cokernel is 1-dimensional. 
Consider $\phi: W \to A^*$ where $W$ is
an indecomposable direct summand of $M$. 
If $W = A$, then $\phi$ factors. 
Suppose $W$ is serial.
Then one easily calculates dimensions and gets that $\phi$ factors. 
If $W = A^*$ and $\phi$ is not an isomorphism, then $\phi$ factors
through the socle quotient, and this is a direct sum of serials. 
Hence $\phi$ factors by what we have already seen.
This shows $\pdim(L(A^*)) \leq 2$.

Next, assume that $X$ is serial but not simple.
Similarly as above, one can show that $\pdim(L(X)) \leq 2$ in this
case.
This uses the short exact sequences 
\[
0 \to M(ba) \to A \to M(aba) \to 0,
\]
\[
0 \to M(ab) \to M(b) \oplus M(aba) \to M(ba) \to 0,
\]
\[
0 \to M(b) \to \field \oplus M(ba) \to M(a) \to 0
\]
with terms in $\add(M)$. 
 
Now consider the projective dimension of $L(\field)$.
We start with the exact sequence 
\[
0 \to D \to M(a) \oplus M(b) \to \field \to 0
\]
of $A$-modules, where $D = A/J_A^2$. 
Applying $\Hom_A(M,-)$ gives the exact sequence 
\[
0 \to \Hom_A(M,D) \to P(M(a)) \oplus P(M(b)) \to P(\field),
\]
which has a 1-dimensional cokernel, namely $L(\field)$.
The exact sequence
\[
0 \to J_A \to A \oplus \field \oplus \field \to D  \to 0
\]
gives rise to the projective resolution
\[
0 \to \Hom_A(M, J_A) \to P(A) \oplus P(\field)
\oplus P(\field) \to \Hom_A(M, D) \to 0.
\]
Hence $\pdim(L(\field)) \leq 3$.
Thus, we get $\gldim(\Gamma) \leq 3$, and therefore $\rdim(A) \leq 3$.
Since the algebra $A$ has infinite representation type, we get
$\rdim(A) \geq 3$ by Auslander's theorem.

The same method works for arbitrary string algebras.


\end{document}